\documentclass[12pt]{amsart}
\usepackage{graphicx}
\usepackage{amsmath}
\usepackage{amsfonts}
\usepackage{amssymb}
\newtheorem{theorem}{Theorem}

\newtheorem{corollary}[theorem]{Corollary}
\newtheorem{definition}[theorem]{Definition}

\newtheorem{remark}[theorem]{Remark}
\newtheorem{lemma}[theorem]{Lemma}

\begin{document}

\title[Overdeterminedness of Functional Equations]%
{The Overdeterminedness of a Class of Functional Equations}
\author{Orr Moshe Shalit}
\thanks{{\it key words and phrases.} functional equations on restricted domains, overdeterminedness.}

\maketitle

\begin{abstract}
{We prove a uniqueness theorem for a large class of functional
equations in the plane, which resembles in form a classical result
of Acz$\acute{e}$l. It is also shown that functional equations in
this class are \emph{overdetermined} in the sense of Paneah. This
means that the solutions, if they exist, are determined by the
corresponding relation being fulfilled not in the original domain
of validity, but only at the points of a subset of the boundary of
the domain of validity.}
\end{abstract}

\section{Introduction}

Recall the classical Cauchy functional equation:
\begin{equation}\label{eq:ClassicalCauchyEq}
f(x + y) = f(x) + f(y) \, \, .
\end{equation}
To \emph{solve} the
functional equation usually means : given a set $A \subseteq
\mathbb{R}^2$ and a class of functions $\mathcal{A}$, to find the
family of functions $\mathcal{F} \subseteq \mathcal{A}$ which
consists of all $f$ such that $f(x+y) = f(x) + f(y)$ for all
$(x,y) \in A$. Following Kuczma (\cite{KuczmaRestrictedDomains})
let us call $A$ the \emph{domain of validity}. For example, when
Cauchy first treated (\ref{eq:ClassicalCauchyEq}), he took
$\mathcal{A} = C(\mathbb{R})$, and showed that if the domain of
validity is taken to be $\mathbb{R}^2$ then the set of solutions
to (\ref{eq:ClassicalCauchyEq}) is $\mathcal{F} = \{f : \exists
\lambda . \forall z . f(z) = \lambda z\}$. It has been shown in
various works (see \cite{Sablik1}, \cite{Sablik2}, \cite{Aczel2},
\cite{Pa04OD} and the references therein) that when some
additional smoothness assumptions are imposed on $f$ then even if
the domain of validity is quite small - the graph of an
appropriate function, for example - the set of solutions does not
grow. Thus, using the terminology of Paneah (\cite{Pa04OD}), we
may say that the equation
\begin{equation*}
f(x + y) = f(x) + f(y) \hspace{5 mm} , \hspace{5 mm} (x,y) \in
\mathbb{R}^2
\end{equation*}
is \emph{overdetermined} (for the class of functions satisfying
these additional smoothness assumptions). In
fact, in \cite{Thesis} we proved that the Cauchy equation in
$\mathbb{R}^n$
\begin{displaymath}
f(x_1+y_1, \ldots, x_n+y_n) = f(x_1, \ldots, x_n) + f(y_1, \ldots,
y_n)
\end{displaymath}
is overdetermined for the class $C^1(\mathbb{R}^n,\mathbb{R})$.

One is led to the following questions: (a) \emph{given a class of
functions $\mathcal{A}$, what is the ``smallest'' domain of
validity for which the solutions to (\ref{eq:ClassicalCauchyEq})
are only $f(z) = \lambda z$}, and : (b) \emph{given a domain of
validity, for what $\mathcal{A}$ does the set of solutions to
(\ref{eq:ClassicalCauchyEq}) remain $f(z) = \lambda z$?}

The above questions may be asked with regards to any functional
equation, and it is interesting in general to study how, given a
functional equation, the set of solutions changes when the domain
of validity and the class of functions considered are changed.
For most classical functional equations in 2 variables, the domain
of validity is usually taken to be some large, open set in
$\mathbb{R}^2$. In \cite{Pa04OD} Paneah proved for a sample of
classical functional equations that, under some smoothness
assumptions, their solution is already determined by the
functional equation holding on a much smaller domain of validity,
e.g., a one-dimensional sub-manifold in $\mathbb{R}^2$, and such
equations were called \emph{overdetermined}. In this paper we
consider a class of functional equations that turn out to be
overdetermined (section 3). This class of equations contains some well known
equations such as Jensen's equation and the equation of the
logarithmic mean. Our main
theorem resembles in form (and in fact, was inspired by) a classical
result of Acz$\acute{e}$l. Our proof depends on a minimality result in topological
dynamics to which the next section is devoted.

\section{Some preliminaries in topological dynamics}
A \emph{dynamical system} is a pair $(X, \delta)$, where $X$ is a
metric space and $\delta =
(\delta_1 , \ldots , \delta_N)$ is a set of continuous maps
$\delta_i : X \rightarrow X$. The maps in $\delta$ generate (by
composition) a semigroup of maps $\Phi_\delta$ in the following
manner:
\begin{equation*}
\Phi_\delta^0 = \{\mathbf{id_X}\} \, \, ,
\end{equation*}
\begin{equation*}
\Phi_\delta^m = \{\sigma_1 \circ \cdots \circ \sigma_m \vert
\sigma_1 , \ldots , \sigma_m \in \delta \} \, \, ,
\end{equation*}
and
\begin{equation*}
\Phi_\delta = \bigcup_{m=0}^\infty \Phi_\delta^m \,\,.
\end{equation*}

\begin{definition}
For any $x \in X$, the \emph{orbit} of $x$ is the set
$$O(x) = \{\sigma(x) \vert \sigma \in \Phi_\delta \} \, \, .$$
\end{definition}
\begin{definition}
A dynamical system $(X,\delta)$ is said to be \emph{minimal} if for all $x \in X$
$$\overline{O(x)} = X \, \, .$$
\end{definition}

\begin{lemma} \label{lemma:contracting}
Let $(X,d)$ be a compact, metric space, and let $\delta =
(\delta_1, \delta_2, \ldots , \delta_N )$ be a finite family of
functions $X \rightarrow X$ satisfying
\begin{equation}\label{eq:rangeCoverX}
\delta_1(X) \cup \delta_2(X) \cup \ldots \cup \delta_N(X) = X
\quad .
\end{equation}
If $\delta$ has the property that for all $ i = 1, \ldots , N$ and
all $x,y \in X$
\begin{equation}\label{eq:contracting}
x \neq y \ \Rightarrow \ d(\delta_i(x),\delta_i(y)) < d(x,y)
\end{equation}
then the dynamical system $(X,\delta)$ is minimal.
\end{lemma}

\begin{proof}
Fix $x_0 \in X$ . We must show that for any $y$ in $X$ and
$\epsilon
> 0 $, there is a $z \in O(x_0)$ such that $d(z,y) \leq
\epsilon$. Fix some $y \in X$ and $\epsilon > 0$. The set $S \equiv \{(x_1, x_2) \in X \times X \, \vert \, d(x_1,x_2) \geq \epsilon\}$ is a compact subset of $X \times X$, thus for every $i = 1,2, \ldots , N$, the continuous funtion
$g_i \colon S \to \mathbb{R}$ defined by:
\begin{equation*}
g_i(x_1,x_2) = {d(\delta_i(x_1),\delta_i(x_2)) \over d(x_1,x_2)} \quad , \quad (x_1,x_2) \in S
\end{equation*}
attains a maximum $c_{\epsilon , i}$. By
(\ref{eq:contracting}), \ $c_{\epsilon , i}<1$, for all $i$. Set
$c_\epsilon$ to be the maximum of these constants.

Now choose some $n$
satisfying $c_{\epsilon}^n \cdot {\rm diam}(X) <\epsilon$. Then for all
$\sigma \in \Phi_\delta^n$ and all $x_1, x_2 \in X$
\begin{displaymath}
d(\sigma(x_1),\sigma(x_2)) \leq \epsilon
\end{displaymath}
and thus for all $\sigma \in \Phi_\delta^n$ :
\begin{equation} \label{eq:diam}
{\rm diam}(\sigma(X)) \leq \epsilon \, \, .
\end{equation}

But note that by virtue of (\ref{eq:rangeCoverX}),
\begin{displaymath}
\bigcup_{\sigma \in \Phi_\delta^n} \sigma(X) = X
\end{displaymath}
so that there is an $f \in \Phi_\delta^n$ s.t. $y \in f(X)$. Now
by (\ref{eq:diam}) it follows that for all $x$ it is true that
$d(f(x),y) \leq \epsilon$ so we can choose $z = f(x_0)$ and the
proof is complete.
\end{proof}

\section{The main result}

In 1964 Acz\'{e}l (see \cite{AczelBooklet}) proved the following
uniqueness theorem for a rather wide class of functional
equations:

\begin{theorem}\label{thm:aczel}
Let $f_1, f_2 : I \rightarrow \mathbb{R} $ be continuous solutions
of the equation

\begin{equation}\label{eq:feaczel}
f(F(x,y)) = H[f(x),f(y),x,y] \verb"  ,  " (x,y) \in I^2
\end{equation}
where $I$ is an (open, closed, half-open, finite or infinite)
interval. Suppose that $F:I^2 \rightarrow I$ is continuous and
internal that is,

$$ min(x,y) < F(x,y) < max(x,y) \, \, \emph{\textrm{ if }} \, \,  x \neq y $$
and that either $u \mapsto H(u,v,x,y)$ or $v \mapsto H(u,v,x,y)$
are injections. Further, let $a,b \in I$ and

\begin{displaymath}
f_1(a) = f_2(a) \, \, \emph{\textrm{ and }} \, \, f_1(b) = f_2(b)
\, .
\end{displaymath}
Then
\begin{displaymath}
\forall x \in I . f_1(x) = f_2(x) \, \, .
\end{displaymath}
\end{theorem}

This theorem motivated much work on uniqueness theorems and has
been improved several times. Theorems in the same spirit were
proved for different classes of $F$ and $H$ and for more general
spaces ($\mathbb{R}^2$,$\mathbb{R}^n$, topological vector spaces,
\ldots \footnote{\cite{Aczel2} contains references to these
developments.}). In this section we will prove a refinement of the
above theorem which serves at once both as a uniqueness theorem
for (\ref{eq:feaczel}) and as a proof that all of the equations
that belong to the class treated below are overdetermined.

\begin{theorem}\label{thm:uniqueOD}
Let $I = [a,b]$, $H : \mathbb{R} \times \mathbb{R} \times I \times
I \rightarrow \mathbb{R}$ any function and $F : I^2 \rightarrow I$
a continuous function that satisfies
\begin{description}
  \item[1] $\forall x \neq y . |F(x,b)-F(y,b)|,|F(a,x)-F(a,y)|<|x-y|$
  \item[2] $\exists x_0, y_0 . F(a,x_0) = a \quad {\rm and} \quad F(y_0,b) = b$
\end{description}
For any real $A$ and $B$ there exists at most one solution $f$ to
(\ref{eq:feaczel}) that satisfies the boundary conditions

\begin{equation}\label{eq:boundcond}
f(a) = A \quad , \quad f(b) = B .
\end{equation}
Moreover, if a function $f$ is a solution to (\ref{eq:feaczel})
satisfying (\ref{eq:boundcond}), then it is already determined by
the functional equation

\begin{equation}\label{eq:feongamma}
f(F(x,y)) = H[f(x),f(y),x,y] \quad , \quad (x,y) \in \Gamma
\end{equation}
 where $ \Gamma = ([a,b] \times \{b\} ) \cup (\{a\} \times
[a,b]) $.
\end{theorem}
\begin{proof}
Let us define two maps $\delta_1 , \delta_2 : I \rightarrow I$ by the
formulas
$$\delta_1(x) = F(a,x) \, \, ,$$
$$\delta_2(x) = F(x,b) \, \, .$$
We consider the dynamical system $(I,\delta_1, \delta_2)$.
By the definitions of $\delta_1, \delta_2$ and by the conditions on $F$
we have that
$$\delta_1(b) = \delta_2(a) \, \, ,$$
and that
$$\delta_1(x_0) = a \, \, \textrm{ and } \, \, \delta_2(y_0) = b \, \, ,$$
thus,
$$ \delta_1(I) \cup \delta_2(I) = I . $$
In addition,
$$\forall x \neq y .
|\delta_2(x)-\delta_2(y)|,|\delta_1(x)-\delta_1(y)|<|x-y| \, \, .$$
By lemma \ref{lemma:contracting} it follows that the orbit of any point
in $I$ is dense in $I$.

Now let $f_1$ and $f_2$ be continuous and satisfy
(\ref{eq:boundcond}) and (\ref{eq:feongamma}). We shall show that
for any $z$ in the orbit of $a$
$$ f_1(z) = f_2(z) .$$
For $a$ we already have by (\ref{eq:boundcond}) that
$$f_1(a) = A = f_2(a) .$$
If $z$ is a point for which we know that $f_1(z) = f_2(z)$ then
\begin{displaymath}
f_1(\delta_1(z)) = f_1(F(a,z)) = H[f_1(a),f_1(z),a,z]
\end{displaymath}
by (\ref{eq:feongamma}). But by our assumption on $z$ we can
replace $H[f_1(a),f_1(z),a,z]$ by $H[f_2(a),f_2(z),a,z]$ and
obtain
\begin{displaymath}
f_1(\delta_1(z)) = H[f_2(a),f_2(z),a,z] = f_2(\delta_1(z))
\end{displaymath}
where the last equality follows again from (\ref{eq:feongamma}).
So we have
$$f_1(\delta_1(z)) = f_2(\delta_1(z)) \,\, .$$
Arguing in just the same manner we arrive at the relation
$$f_1(\delta_2(z)) = f_2(\delta_2(z)) \,\,.$$
So all the points in the orbit of $a$ inherit from $a$ the
property of being given the same values by $f_1, f_2$, and so
indeed for any $z \in O (a)$ we have $f_1(z) = f_2(z)$. The
continuity of $f_1 , f_2$ and the density of $O (a)$ imply $f_1 =
f_2$ on $I$.
\end{proof}

\begin{remark}\emph{Note that the above proof suggests an
algorithm that can compute numerically a solution (when such
exists) to a given functional equation on an interval with
boundary data.}
\end{remark}

\begin{remark}\emph{ Note that it follows from the above theorem that
usually (\ref{eq:feaczel}) will not have a solution, even if
(\ref{eq:feongamma}) has a solution.}
\end{remark}

As a corollary of the above theorem we have the overdeterminedness
of Jensen's functional equation.
\begin{corollary}
Let $\alpha$ and $\beta$ be two positive numbers satisfying
$\alpha+\beta = 1$, and let $I=[a,b]$ be some closed interval.
Then all continuous solutions $f$ of the functional equation
\begin{displaymath}
f(\alpha x + \beta y) = \alpha f(x) + \beta f(y) \,\, , \,\, (x,y)
\in I^2
\end{displaymath}
are of the form
\begin{displaymath}
f(z) = \lambda z + \mu
\end{displaymath}
for some constants $\lambda, \mu \in \mathbb{R}$. Moreover, these
solutions are already determined by the functional equation
\begin{displaymath}
f(\alpha x + \beta y) = \alpha f(x) + \beta f(y) \,\, , \,\, (x,y)
\in \Gamma
\end{displaymath}
where $ \Gamma = ([a,b] \times \{b\} ) \cup (\{a\} \times [a,b])
$.
\end{corollary}

\section{Acknowledgments} The results in this paper are part of the
author's MSc. thesis which was done under the supervision of Boris Paneah,
to whom I am grateful for introducing me to the
problem and for advising and encouraging me along the way.
Yoav Benyamini and Daniel Reem suggested various improvements to the text.

\bibliographystyle{abbrv}
\bibliography{OD_theorem}

\vspace{15 pt} {\noindent Orr Moshe Shalit, Department of
Mathematics, Technion\\32000 Haifa, Israel. e-mail:
orrms@tx.technion.ac.il}

\end{document}